\documentclass[11pt,a4paper]{amsart}
\usepackage{amssymb}
\usepackage{amsmath}
\usepackage{amsfonts}
\usepackage{bbm}
\usepackage{enumerate}

\newtheorem{theorem}{Theorem}[section]

\theoremstyle{definition}\newtheorem*{proof0}{Proof}

\theoremstyle{definition}

\theoremstyle{definition}

\begin{document}

\title{Identities of Eulerain and ordered Stirling numbers over a multiset}

\author{Joon Yop Lee}
\address{Department of Mathematics, POSTECH, Pohang 790-784, Korea}
\email{flutelee@postech.ac.kr}

\begin{abstract}
By considering Eulerian numbers and ordered Stirling numbers of the second and third kinds over a multiset, we generalize identities of Eulerian numbers and Stirling numbers of the second and third kinds and provide $q$-analogs of these generalizations. Using these generalizations, we also compute Eulerian numbers and ordered Stirling numbers of the second kind over a multiset.

\end{abstract}

\maketitle

\section{Introduction} 
Let $\left\langle {d\atop i}\right\rangle$
be an Eulerian number, which is the number of permutations of $\{1,2,\ldots,d\}$ with $i$ descents. A well-known identity of Eulerian numbers is Worpitzky's identity \cite{worpitzky}
\begin{equation}\label{identity1.1}
x^d=\sum_{i=0}^{d-1}\left\langle {d\atop i}\right\rangle{x-1-i+d\choose d},
\end{equation}
and the following Carlitz's identity \cite{carlitz1,carlitz2} is a $q$-analog of this identity
\begin{equation}\label{identity1.11}
\bigg[{x\atop 1}\bigg]_q^d=\sum_{i=1}^dA_{d,i}(q)\bigg[{x-1+i\atop d}\bigg]_q
\end{equation}
where $\big[{x\atop m}]_q=\prod_{i=1}^m(1-q^{x-m+i})/(1-q^i)$ and $A_{d,i}(q)$ is a polynomial of $q$.
Eulerian numbers are defined over an ordinary set, thus one way to generalize identities (\ref{identity1.1}) and (\ref{identity1.11}) is to consider Eulerian numbers over a multiset. By considering Eulerian numbers over a multiset, we may be able to generalize identities (\ref{identity1.1}) and (\ref{identity1.11}). 

We can apply such an idea to identities of Stirling numbers of the second and third kinds.
Let ${d\brace k}$ be a Stirling number of the second kind, which is the number of partitions of a $d$ element set into $k$ nonempty sets. String numbers of the second kind satisfy the following polynomial identity
\begin{equation}\label{identity1.2}
x^d=\sum_{k=1}^d k!\left\{{d\atop k}\right\}{x\choose k}
\end{equation}
A Stirling number of the third kind \cite{tsylova} (or a Lah number), denoted by $\left\lfloor{d\atop k}\right\rfloor$, is the number of ways to partition a $d$ element set into $k$ nonempty linearly ordered subsets \cite{aigner, motzkin}. For example, the partitions of $\{1,2,3\}$ into $2$ nonempty linearly ordered sets are $\{(1), (2, 3)\}$, $\{(1), (3, 2)\}$, $\{(2), (1, 3)\}$, $\{(2), (3, 1)\}$, $\{(3), (1, 2)\}$, and $\{(3), (2, 1)\}$.
Lah \cite{lah} defined Stirling numbers of the third kind by the following polynomial identity
\begin{equation}\label{identity1.3}
d!{x-1+d\choose d}=\sum_{k=1}^d k!\left\lfloor {d\atop k}\right\rfloor{x\choose k}.
\end{equation}
By considering Stirling numbers of the second and third kinds over a multiset, we may also be able to generalize identities (\ref{identity1.2}) and (\ref{identity1.3}) and provide $q$-analogs of these generalizations.

Motivated by these ideas, we will generalize identities (\ref{identity1.1}), (\ref{identity1.2}), and (\ref{identity1.3}) and provide $q$-analogs of these generalizations. In particular, we generalize identities (\ref{identity1.2}) and (\ref{identity1.3}) and obtain $q$-analogs of these generalizations by considering ordered Stirling numbers of the second kind $k!{d\brace k}$ and third kind $k!\left\lfloor{d\atop k}\right\rfloor$ over a multiset, rather than ordinary Stirling numbers of the second and third kinds.
Using generalizations of identities (\ref{identity1.1}) and (\ref{identity1.2}), we also compute Eulerian numbers and ordered Stirling numbers of the second kind over a  multiset.

The gist of our idea is the following. For a sequence of finite sets of lattice points $S_0$, $S_1$, $S_2$,\,\ldots, we compute the numbers of elements in $S_n$ by two different ways and obtain a polynomial identity.
To obtain a $q$-analog of this identity, we compute the following generating function
$$\sum_{(x_1,x_2,\ldots,x_d)\in S_n}
q^{x_1+x_2+\cdots+x_d}$$
by two different ways.

\section{A triangulation of the product of simplexes}

In this section, we introduce a triangulation of the product of simplexes. This triangulation will be the main tool to generalize identities (\ref{identity1.1}), (\ref{identity1.2}), and (\ref{identity1.3}) and provide $q$-analogs of these generalizations.

For a positive integer $n$, let $[n]=\{1,2,\ldots,n\}$ and $[n]_0=\{0,1,\ldots,n\}$. We denote a point in $\mathbb{R}^d$ by $\mathbf{x}=(x_1,x_2,\ldots,x_d)$. We denote the zero vector of $\mathbb{R}^d$ by $\mathbf{e}_{d,0}$ and for each $i\in[d]$ we denote the $i$th unit vector of $\mathbb{R}^d$ by $\mathbf{e}_{d,i}=(0,\ldots,0,1,0,\ldots,0)$ (with $1$ in the $i$th coordinate).  
 We define $\alpha^d_0$ to be a $d$-dimensional simplex whose vertexes are $\mathbf{e}_{d,0}$, $\mathbf{e}_{d,0}+\mathbf{e}_{d,1}$,\,\ldots,  $\mathbf{e}_{d,0}+\mathbf{e}_{d,1}+\cdots+\mathbf{e}_{d,d}$. Note that $\alpha^d_0$ is the set of points $\mathbf{x}$ such that $1\ge x_1\ge x_2\ge\cdots\ge x_d\ge 0$.

Let $d=d_1+d_2+\cdots+d_l$ be a sum of nonnegative integers.
For each $j\in[l]$ writing a point of $\mathbb{R}^{d_j}$ by $\mathbf{x}_j=(x_{j,1},x_{j,2},\ldots,x_{j,d_j})$, we also denote a point of $\mathbb{R}^d=\prod_{j=1}^l\mathbb{R}^{d_j}$ by $\mathbf{x}=\prod_{j=1}^l\mathbf{x}_j$. 
Letting $\{(j_h,i_h)\,|\,h\in[d]\}=\bigcup_{j=1}^l\{j\}\times[d_j]$, for a point $\mathbf{x}$ in $\mathbb{R}^d$ we always assume that
\begin{enumerate}[1.]
\item
$ x_{j_1,i_1}\ge x_{j_2,i_2}\ge\cdots\ge x_{j_d,i_d}$,

\item if $x_{j_h,i_h}=x_{j_{h+1},i_{h+1}}$ then either $j_h<j_{h+1}$ or $(j_h,i_h)=(j_{h+1},i_{h+1}-1)$.
\end{enumerate}
Let $\mathbf{d}=(d_1,d_2,\ldots,d_l)$. We denote the product of simplexes $\alpha^{d_1}_0$, $\alpha^{d_2}_0$,\,\ldots, $\alpha^{d_l}_0$ by $\alpha^{\mathbf{d}}=\prod_{j=1}^l\alpha^{d_j}_0$ and the vertexes of $\alpha^{\mathbf{d}}$ by $\mathbf{e}_{i_1,i_2,\ldots,i_l}=\prod_{j=1}^l(\mathbf{e}_{d_j,0}+\mathbf{e}_{d_j,1}+\cdots+\mathbf{e}_{d_j,i_j})$ where $(i_1,i_2,\ldots,i_l)$ ranges over the set $\prod_{j=1}^l[d_j]_0$.

Let $S(\mathbf{d})=S(d_1,d_2,\ldots,d_l)$ be a multiset such that for each $j\in[l]$ the number of $j$ in $S(\mathbf{d})$ is $d_j$. We define $\mathfrak{S}(\mathbf{d})$ to be the permutation set of $S(\mathbf{d})$ and denote a permutation of $S(\mathbf{d})$ by $\sigma=[\sigma_1,\sigma_2,\ldots,\sigma_d]$.
For a permutation $\sigma$ of $S(\mathbf{d})$ we define $\alpha^d(\sigma)$ to be a $d$-simplex with the vertexes $\mathbf{e}_{i_{h,1},i_{h,2},\ldots,i_{h,l}}$ for $h\in[d]_0$ such that $(i_{0,1},i_{0,2},\ldots,i_{0,l})=\mathbf{e}_{l,0}$ and $(i_{h,1},i_{h,2},\ldots,i_{h,l})=\sum_{a=1}^h\mathbf{e}_{l,\sigma_a}$ for $h\in[d]$.
By definition, $\alpha^d(\sigma)$ is the set of points $\mathbf{x}$ such that
\begin{enumerate}[1.]
\item
$ x_{\sigma_1,i_1}\ge x_{\sigma_2,i_2}\ge\cdots\ge x_{\sigma_d,i_d}$,

\item if $x_{\sigma_h,i_h}=x_{\sigma_{h+1},i_{h+1}}$ then either $\sigma_h<\sigma_{h+1}$ or $(\sigma_h,i_h)=(\sigma_{h+1},i_{h+1}-1)$.
\end{enumerate}

Let $\mathbf{x}$ be a point of $\alpha^{\mathbf{d}}$.
Letting $x_{j_0,i_0}=1$ and $x_{j_{d+1},i_{d+1}}=0$, we can represent  $\mathbf{x}$ by the following sum
$$\mathbf{x}=\sum_{k=0}^d(x_{j_h,i_h}-x_{j_{h+1},i_{h+1}})\mathbf{e}_{i_{h,1},i_{h,2},\ldots,i_{h,l}}$$
where $(i_{0,1},i_{0,2},\ldots,i_{0,l})=\mathbf{e}_{l.0}$ and $(i_{h,1},i_{h,2},\ldots,i_{h,l})=\sum_{a=1}^h\mathbf{e}_{l,j_a}$ for $h\in[d]$.
This implies that $\mathbf{x}$ is a convex sum of the vertexes of $\alpha^d([j_1,j_2,\ldots,j_d])$. Since this is true for every point of $\alpha^{\mathbf{d}}$, it follows that $\alpha^{\mathbf{d}}$ is the union of the $d$-simplexes $\alpha^d(\sigma)$ for $\sigma\in\mathfrak{S}(\mathbf{d})$.
Moreover, the set composed of the $d$-simplexes $\alpha^d(\sigma)$ for $\sigma\in\mathfrak{S}(\mathbf{d})$ and their faces is a triangulation of $\alpha^{\mathbf{d}}$ by definition. We denote this triangulation by $\mathcal{T}_{\alpha^{\mathbf{d}}}$.

\section{Eulerian numbers over a multiset}

Let $R$ be a subset of $\mathbb{R}^d$. For a nonnegative integer $n$ we denote $nR=\{nr\,|\,r\in R\}$ and we define $\mathbb{Z}(R)$ to be the set of lattice points in $R$. Recall that a lattice point is a point with nonnegative integers entries. For a point $\mathbf{x}$ let $q^{\mathbf{x}}=q^{x_1+x_2+\cdots+x_d}$. To generalize identity (\ref{identity1.1}), we will compute the number of elements in  $\mathbb{Z}(n\alpha^{\mathbf{d}})$ by two different ways and to obtain a $q$-analog of this generalization, which is a generalization of identity (\ref{identity1.11}), we will compute $f_1(q)=\sum_{\mathbf{x}\in \mathbb{Z}(n\alpha^{\mathbf{d}})}q^{\mathbf{x}}$ by two different ways.
 
For a permutation $\sigma$ of $S(\mathbf{d})$ let $D(\sigma)$ be the descent set of $\sigma$, that is, the set of indexes $h$ such that $\sigma_h>\sigma_{h+1}$.
We define $A(\sigma)$ to be a subset of $\alpha^d(\sigma)$ composed of points $\mathbf{x}$ such that $x_{\sigma_h,i_h}>x_{\sigma_{h+1},i_{h+1}}$ for all $h$ where $h$ ranges over $D(\sigma)$. Then for each point $\mathbf{x}$ in $\alpha^{\mathbf{d}}$ there is a unique permutation $\sigma$ such that $A(\sigma)$ contains $\mathbf{x}$. Therefore $\alpha^{\mathbf{d}}$ is the disjoint union of the sets $A(\sigma)$ for $\sigma\in\mathfrak{S}(\mathbf{d})$. We call this disjoint union of $\alpha^{\mathbf{d}}$ the \emph{first decomposition} of $\alpha^{\mathbf{d}}$. Stanley \cite{stanley} obtained such a decomposition when $d_1=d_2=\cdots=d_l=1$, that is, $\alpha^{\mathbf{d}}$ is an $l$-dimensional hypercube.

The set $\mathbb{Z}(n\alpha^{\mathbf{d}})$ is the product of $\mathbb{Z}(n\alpha^{d_1})$, $\mathbb{Z}(n\alpha^{d_2})$,\,\ldots, $\mathbb{Z}(n\alpha^{d_l})$ and for each $j\in[l]$ the number of elements in the set $\mathbb{Z}(n\alpha^{d_j})$ is ${n+d_j\choose d_j}$, thus
$$|\mathbb{Z}(n\alpha^{\mathbf{d}})|=
\prod_{j=1}^l|\mathbb{Z}(n\alpha^{d_j})|=\prod_{j=1}^l{n+d_j\choose d_j}.$$
A point $\mathbf{x}$ is an element of $\mathbb{Z}(nA(\sigma))$ if and only $\mathbf{x}$ is a lattice point such that
\begin{enumerate}[1.]
\item $n\ge x_{\sigma_1,i_1}\ge x_{\sigma_2,i_2}\ge\cdots\ge x_{\sigma_d,i_d}\ge 0$,

\item $x_{\sigma_h,i_h}\ge x_{\sigma_{h+1},i_{h+1}}+1$ whenever $h$ is a descent of $\sigma$.
\end{enumerate}
Thus the number of elements in $\mathbb{Z}(nA(\sigma))$ is ${n-|D(\sigma)|+d\choose d}$. Therefore if we denote  by $\left\langle {\mathbf{d}\atop i}\right\rangle$ the number of permutations of $S(\mathbf{d})$ with $i$ descents, then by the first decomposition of $\alpha^{\mathbf{d}}$ we obtain
\begin{align*}
|\mathbb{Z}(n\alpha^{\mathbf{d}})|
&=|\mathbb{Z}(\biguplus_{\sigma\in\mathfrak{S}(\mathbf{d})}nA(\sigma))|=\sum_{\sigma\in\mathfrak{S}(\mathbf{d})}|\mathbb{Z}(nA(\sigma))|\nonumber\\
&=\sum_{\sigma\in\mathfrak{S}(\mathbf{d})}{n-|D(\sigma)|+d\choose d}=\sum_{i=1}^{d-1}\left\langle{\mathbf{d}\atop i}\right\rangle
{n-i+d\choose d}\nonumber.\end{align*}
As a result, we obtain the following identity
\begin{equation}\label{identity3.1}
\prod_{j=1}^l{n+d_j\choose d_j}=\sum_{i=0}^{d-1}\left\langle{\mathbf{d}\atop i}\right\rangle
{n-i+d\choose d}.
\end{equation}
Note that Kim and Lee [Polytope numbers and their properties, arXiv:1206.0511] derived identity (\ref{identity3.1}) by using the concept of polytope numbers.
Since identity (\ref{identity3.1}) is true for every nonnegative integer $n$, it follows that
\begin{equation}\label{identity3.2}
\prod_{j=1}^l{x+d_j\choose d_j}=\sum_{i=0}^{d-1}\left\langle{\mathbf{d}\atop i}\right\rangle
{x-i+d\choose d}.
\end{equation}
\begin{theorem}
Eulerain numbers over a multiset $S(\mathbf{d})$ satisfies
$$\prod_{j=1}^l{x+d_j\choose d_j}=\sum_{i=0}^{d-1}\left\langle{\mathbf{d}\atop i}\right\rangle
{x-i+d\choose d}.$$
\end{theorem}
When $d_1=d_2=\cdots=d_l=1$, identity (\ref{identity3.2}) becomes identity (\ref{identity1.1}). Thus identity (\ref{identity3.2}) is a generalization of Worpitzky's identity \cite{worpitzky}.

Using identity (\ref{identity3.2}), we can obtain Eulerian numbers over $S(\mathbf{d})$ in the following way.
Let $M_1$ be a $d\times d$ matrix defined by
$$M_1=\begin{bmatrix}
{d\choose d}&0&\cdots&0\\
{d+1\choose d}&{d\choose d}&\cdots&0\\
\cdots&\cdots&\cdots&\cdots\\
{2d-1\choose d}&{2d-2\choose d}&\cdots&{d\choose d}
\end{bmatrix}.$$
Then from identity (\ref{identity3.1}) we can obtain the following 
matrix identity
$$\begin{bmatrix}
\prod_{j=1}^l{d_j\choose d_j}\\
\prod_{j=1}^l{1+d_j\choose d_j}\\
\cdots\\
\prod_{j=1}^l{d-1+d_j\choose d_j}
\end{bmatrix}=
\begin{bmatrix}
{d\choose d}&0&\cdots&0\\
{d+1\choose d}&{d\choose d}&\cdots&0\\
\cdots&\cdots&\cdots&\cdots\\
{2d-1\choose d}&{2d-2\choose d}&\cdots&{d\choose d}
\end{bmatrix}
\begin{bmatrix}
\Big\langle{\mathbf{d}\atop 0}\Big\rangle\\
\left\langle{\mathbf{d}\atop 1}\right\rangle\\
\cdots\\
\left\langle{\mathbf{d}\atop d-1}\right\rangle
\end{bmatrix}.
$$
Since the inverse of $M_1$ is
$$M_1^{-1}=
\begin{bmatrix}
(-1)^0{d+1\choose 0}&0&\cdots&0\\
(-1)^1{d+1\choose 1}&(-1)^0{d+1\choose 0}&\cdots&0\\
\cdots&\cdots&\cdots&\cdots\\
(-1)^{d-1}{d+1\choose d-1}&(-1)^{d-2}{d+1\choose d-2}&\cdots&(-1)^0{d+1\choose 0}
\end{bmatrix}$$
by the Gaussian elimination,
it follows that
$$\bigg\langle{\mathbf{d}\atop i}\bigg
\rangle
=\sum_{h=0}^i(-1)^{i-h}{d+1\choose i-h}
\prod_{j=1}^l{h+d_j\choose d_j}.$$
\begin{theorem}
Eulerian numbers over a multiset $S(\mathbf{d})$ are given by
$$\bigg\langle{\mathbf{d}\atop i}\bigg
\rangle
=\sum_{h=0}^i(-1)^{i-h}{d+1\choose i-h}
\prod_{j=1}^l{h+d_j\choose d_j}.$$
\end{theorem}

Now we consider $f_1(q)$.
By definition, $\mathbf{x}_j$ is in $\mathbb{Z}(n\alpha^{d_j}_0)$ if and only if $\mathbf{x}_j$ is a lattice point such that $n\ge x_{j,1}\ge x_{j,2}\ge\cdots\ge x_{j,d_j}\ge 0$, thus the generating function $\sum_{\mathbf{x}_j\in\mathbb{Z}(n\alpha^{d_j}_0)}q^{\mathbf{x}_j}$ is a $q$-binomial coefficient $\Big[{n+d_j\atop d_j}\Big]_q$ \cite{andrews,sylvester}.
Since $\mathbb{Z}(n\alpha^{\mathbf{d}})=\prod_{j=1}^l\mathbb{Z}(n\alpha^{d_j}_0)$ and $q^{\mathbf{x}}=\prod_{j=1}^lq^{\mathbf{x}_j}$, it follows that 
\begin{equation}\label{identity3.3}
f_1(q)=\prod_{j=1}^l\sum_{\mathbf{x}_j\in\mathbb{Z}(n\alpha^{d_j}_0)}q^{\mathbf{x}_j}\\
=\prod_{j=1}^l\bigg[{n+d_j\atop d_j}\bigg]_q.
\end{equation}

For a permutation $\sigma$ of $S(\mathbf{d})$ we define the \emph{major index} of $\sigma$ to be $maj(\sigma)=\sum_{j\in D(\sigma)}j$. 
A lattice point $\mathbf{x}$ is in $\mathbb{Z}(n A(\sigma))$ if and only if
$$\begin{cases}
x_{\sigma_h,i_h}\ge x_{\sigma_{h+1},i_{h+1}}&\text{for}\ h\in[d+1]\setminus D(\sigma)\\ x_{\sigma_h,i_h}\ge x_{\sigma_{h+1},i_{h+1}}+1&\text{for}\ h\in D(\sigma)
\end{cases}.$$
Therefore
$$\sum_{\mathbf{x}\in \mathbb{Z}(n A(\sigma))}q^{\mathbf{x}}=\sum_{\mathbf{x}\in \mathbb{Z}(n A(\sigma))}\prod_{k=1}^dq^{x_{j_k,i_k}}=q^{maj(\sigma)}\bigg[{n-|D(\sigma)|+d\atop d}\bigg]_q.$$
Since $\mathbb{Z}(n\alpha^{\mathbf{d}})$ is the disjoint union of the sets $\mathbb{Z}(nA(\sigma))$ for $\sigma\in\mathfrak{S}(\mathbf{d})$, we obtain
\begin{align}\label{identity3.4}
f_1(q)=&\sum_{\sigma\in\mathfrak{S}(\mathbf{d})}
\sum_{\mathbf{x}\in\mathbb{Z}(nA(\sigma))}q^{\mathbf{x}}\nonumber\\
=&\sum_{\sigma\in\mathfrak{S}(\mathbf{d})}q^{maj(\sigma)}\bigg[{n-|D(\sigma)|+d\atop d}\bigg]_q\nonumber\\
=&\sum_{i=1}^dA_{\mathbf{d},i}(q)\bigg[{n+i\atop d}\bigg]_q
\end{align}
where $A_{\mathbf{d},i}(q)=\sum_{\sigma\in\mathfrak{S}(\mathbf{d})\atop |D(\sigma)|=d-i}q^{maj(\sigma)}$.
Therefore combining identities (\ref{identity3.3}) and (\ref{identity3.4}), we obtain 
\begin{equation}\label{identity3.5}
\prod_{j=1}^l\bigg[{n+d_j\atop d_j}\bigg]_q=\sum_{i=1}^{d}A_{\mathbf{d},i}(q)\bigg[{n+i\atop d}\bigg]_q.\end{equation}
Identity (\ref{identity3.5}) is true for every nonnegative integer $n$, thus by the substitution $n\rightarrow x$ we obtain the 
the following $q$-analog of identity (\ref{identity3.2})
\begin{equation}\label{identity3.6}
\prod_{j=1}^l\bigg[{x+d_j\atop d_j}\bigg]_q=\sum_{i=1}^dA_{\mathbf{d},i}(q)\bigg[{x+i\atop d}\bigg]_q.
\end{equation}
\begin{theorem}
A $q$-analog of identity (\ref{identity3.2}) is
$$\prod_{j=1}^l\bigg[{x+d_j\atop d_j}\bigg]_q=\sum_{i=1}^dA_{\mathbf{d},i}(q)\bigg[{x+i\atop d}\bigg]_q.$$
\end{theorem}

Note that identity (\ref{identity3.6}) is a generalization of Carlitz\rq{}s identity (\ref{identity1.11}).

\section{Ordered Stirling numbers of the second kind over a multiset}

Similar to the case of Eulearin numbers over a multiset, to generalize identity (\ref{identity1.2}), we will compute the number of elements in $\mathbb{Z}(n\alpha^{\mathbf{d}})$ by two different ways  and  to obtain a $q$-analog of this generalization, we will compute $\sum_{\mathbf{x}\in \mathbb{Z}(n\alpha^{\mathbf{d}})}q^{\mathbf{x}}$ by two different ways.

For two vectors $\mathbf{v}=(v_1,v_2,\ldots,v_l)$ and $\mathbf{v}^{\prime}=(v^{\prime}_1,v^{\prime}_2,\ldots,v^{\prime}_l)$, we define $\mathbf{v}\leq \mathbf{v}^{\prime}$ if $v_i\leq v^{\prime}_i$ for all $i\in[l]$. 
Let  $\mathcal{I}(\alpha^{\mathbf{d}})$ be the set of simplexes in $\mathcal{T}_{\alpha^{\mathbf{d}}}$ that contains
$\mathbf{e}_{0,0,\ldots,0}$ and $\mathbf{e}_{d_1,d_2,\ldots,d_l}$.
For a $k$-dimensional simplex $\alpha^k$ ($k\in[d]$) in $\mathcal{I}(\alpha^{\mathbf{d}})$ let $\{\mathbf{e}_{i_{h,1},i_{h,2},\ldots,i_{h,l}}\,|\,h\in[k]_0\}$ be the vertex set of $\alpha^k$ where
\begin{equation}\label{condition4.1}
\begin{cases}
(i_{0,1},i_{0,2},\ldots,i_{0,l})=\mathbf{0}\\
(i_{h,1},i_{h,2},\ldots,i_{h,l})<(i_{h+1,1},i_{h+1,2},\ldots,i_{h+1,l})\ \text{for}\ h\in[k-1]_0\\
(i_{k,1},i_{k,2},\ldots,i_{k,l})=(d_1,d_2,\ldots,d_l).
\end{cases}.
\end{equation}
We define $I(\alpha^k)$ to be the set of convex sums $\mathbf{x}=\sum_{h=0}^kc_h\mathbf{e}_{i_{h,1},i_{h,2},\ldots,i_{h,l}}$ of the vertexes of $\alpha^k$ such that
$c_h>0$ for $h\in[k-1]$.
This is implies that if $\mathbf{x}$ is a point of $I(\alpha^k)$, then there are $k-1$ numbers $h_1$, $h_2$,\,\ldots, $h_{k-1}$ in $[d]$ such that
\begin{equation}\label{interior}
\begin{cases}
x_{j_1,i_1}=\cdots=x_{j_{h_1},i_{h_1}}>\\
x_{j_{h_1+1},i_{h_1+1}}=\cdots=x_{j_{h_2},i_{h_2}}>\\
\hdotsfor{1}\\
x_{j_{h_{k-1}+1},i_{h_{k-1}+1}}=\cdots=x_{j_d,i_d}
\end{cases}.\end{equation}

Let $\mathbf{x}$ be a point of $\alpha^{\mathbf{d}}$.
Suppose that $\mathbf{x}$ satisfies condition (\ref{interior}).
Then a simplex $\alpha^k$ in $\mathcal{I}(\alpha^{\mathbf{d}})$ with the vertex set $\{\mathbf{e}_{i_{a,1},i_{a,2},\ldots,i_{a,l}}\,|\,a\in[k]_0\}$ satisfying
$$\begin{cases}
(i_{0,1},i_{0,2},\ldots,i_{0,l})=\mathbf{0}\\
(i_{a,1},i_{a,2},\ldots,i_{a,l})=\sum_{b=1}^{h_a}\mathbf{e}_{l,j_b}\ \text{for}\ a\in[k-1]\\
(i_{k,1},i_{k,2},\ldots,i_{k,l})=(d_1,d_2,\ldots,d_l)
\end{cases}$$
is the unique simplex such that $I(\alpha^k)$ contains $\mathbf{x}$. Since $I(\alpha^k)$ is a subset of $\alpha^{\mathbf{d}}$, it follows that $\alpha^{\mathbf{d}}$ is the disjoint union of the sets $I(\alpha^k)$ for $\alpha^k\in\mathcal{I}(\alpha^{\mathbf{d}})$. We call this disjoint union of $\alpha^{\mathbf{d}}$ the \emph{second decomposition} of $\alpha^{\mathbf{d}}$

Let $\alpha^k$ be a $k$-dimensional simplex in $\mathcal{I}(\alpha^{\mathbf{d}})$ with the vertex set $\{\mathbf{e}_{i_{h,1},i_{h,2},\ldots,i_{h,l}}\,|\,h\in[k]_0\}$ that satisfies condition (\ref{condition4.1}).
We define the \emph{major index} of $\alpha^k$ to be $maj(\alpha^k)=
\sum_{h=1}^{k-1}\sum_{j=1}^li_{h,j}$.
For each $h\in[k]$ if we denote $S_h=S(i_{h,1},i_{h,2},\ldots,i_{h,l})\setminus S(i_{h-1,1},i_{h-1,2},\ldots,i_{h-1,l})$, then $(S_1,S_2,\ldots,S_k)$ is an ordered partition of $S(\mathbf{d})$ into $k$ nonempty multisets. Therefore
the number of $k$-dimensional simplexes in $\mathcal{I}(\alpha^{\mathbf{d}})$ is the number of ordered partitions of $S(\mathbf{d})$ into $k$ nonempty multisets. We call this number an \emph{ordered Stirling number of the second kind} and denote it by ${\mathbf{d}\brace k}_O$. Note that if $d_1=d_2=\cdots=d_l=1$, then ${\mathbf{d}\brace k}_O=k!{l\brace k}$.

By definition, the number of elements in the set $\mathbb{Z}(nI(\alpha^k))$ is ${n+1\choose k}$, thus by the second decomposition of $\alpha^{\mathbf{d}}$ it follows that
\begin{align*}
\prod_{j=1}^l{n+d_j\choose d_j}&=\prod_{j=1}^l|\mathbb{Z}(n\alpha^{d_j})|=|\mathbb{Z}(n\alpha^{\mathbf{d}})|\\
&=|\mathbb{Z}(n\biguplus_{k=1}^d\biguplus_{\alpha^k\in\mathcal{I}(\alpha^{d_1,\ldots,d_l})}I(\alpha^k))|\\
&=\sum_{k=1}^d\sum_{\alpha^k\in\mathcal{I}(\alpha^{\mathbf{d}})}|\mathbb{Z}(nI(\alpha^k))|\\
&=\sum_{k=1}^d{\mathbf{d}\brace k}_O{n+1\choose k}.
\end{align*}
As a result, by the substitution $n\rightarrow x$ we obtain
\begin{equation}\label{identity4.1}
\prod_{j=1}^l{x+d_j\choose d_j}=\sum_{k=1}^d{\mathbf{d}\brace k}_O{x+1\choose k}.
\end{equation}
\begin{theorem}
Ordered Stirling numbers of the second kind over $S(\mathbf{d})$ satisfies 
$$\prod_{j=1}^l{x+d_j\choose d_j}=\sum_{k=1}^d{\mathbf{d}\brace k}_O{x+1\choose k}.$$
\end{theorem}

Similar to the case of Eulerian numbers over a multiset, we can obtain ordered Stirling numbers of the second kind over a multiset as follows. Let $M_2$ be a $d\times d$ matrix defined by
$$M_2=\begin{bmatrix}
{1\choose 1}&0&\cdots&0\\
{2\choose 1}&{2\choose 2}&\cdots&0\\
\cdots&\cdots&\cdots&\cdots\\
{d\choose 1}&{d\choose 2}&\cdots&{d\choose d}
\end{bmatrix},$$
Then we can rewrite identity (\ref{identity4.1}) in the following form
$$\begin{bmatrix}
\prod_{j=1}^l{d_j\choose d_j}\\
\prod_{j=1}^l{1+d_j\choose d_j}\\
\cdots\\
\prod_{j=1}^l{d-1+d_j\choose d_j}
\end{bmatrix}=
\begin{bmatrix}
{1\choose 1}&0&\cdots&0\\
{2\choose 1}&{2\choose 2}&\cdots&0\\
\cdots&\cdots&\cdots&\cdots\\
{d\choose 1}&{d\choose 2}&\cdots&{d\choose d}
\end{bmatrix}
\begin{bmatrix}
{\mathbf{d}\brace 1}_O\\
{\mathbf{d}\brace 2}_O\\
\cdots\\
{\mathbf{d}\brace d}_O
\end{bmatrix}.$$
Since by the Gaussian elimination the inverse of $M_2$ is
$$M_2^{-1}=
\begin{bmatrix}
(-1)^{0}{1\choose 0}&0&\cdots&0\\
(-1)^{1}{2\choose 1}&(-1)^{0}{2\choose 0}&\cdots&0\\
\cdots&\cdots&\cdots&\cdots\\
(-1)^{d-1}{d\choose d-1}&(-1)^{d-2}{d\choose d-2}&\cdots&(-1)^{0}{d\choose 0}
\end{bmatrix},$$
it follows that
$${\mathbf{d}\brace k}_O
=\sum_{h=0}^{k-1}(-1)^{k-1-h}{k\choose h}
\prod_{j=1}^l{h+d_j\choose d_j}
.$$
\begin{theorem}
Stirling numbers of the second kind over $S(\mathbf{d})$ are given by
$${\mathbf{d}\brace k}_O
=\sum_{h=0}^{k-1}(-1)^{k-1-h}{k\choose h}
\prod_{j=1}^l{h+d_j\choose d_j}.$$
\end{theorem}

For each $\alpha^k\in\mathcal{I}(\alpha^{\mathbf{d}})$ the sum of $q^{\mathbf{x}}$ for $\mathbf{x}\in\mathbb{Z}(nI(\alpha^k))$ is $q^{maj(\alpha^k)}\Big[{n+1\atop k}\Big]_q$. Thus it follows that 
\begin{align*}
f_1(q)
&=\sum_{k=1}^d\sum_{\alpha^k\in\mathcal{I}(\alpha^{\mathbf{d}})}\sum_{\mathbf{x}\in\mathbb{Z}(nI(\alpha^k))}
q^{\mathbf{x}}\\
&=\sum_{k=1}^d\sum_{\alpha^k\in\mathcal{I}(\alpha^{\mathbf{d}})}q^{maj(\alpha^k)}\bigg[{n+1\atop k}\bigg]_q\\
&=\sum_{k=1}^dB_{\mathbf{d},k}(q)\bigg[{n+1\atop k}\bigg]_q
\end{align*}
where $B_{\mathbf{d},k}(q)=\sum_{\alpha^k\in\mathcal{I}(\alpha^{\mathbf{d}})}q^{maj(\alpha^k)}$.
In addition, $f_1(q)=\prod_{j=1}^l\Big[{n+d_j\atop d_j}\Big]_q$, therefore we obtain
$$\prod_{j=1}^l\bigg[{n+d_j\atop d_j}\bigg]_q
=\sum_{k=1}^dB_{\mathbf{d},k}(q)\bigg[{n+1\atop k}\bigg]_q.$$
As a result, by the substitution $n\rightarrow x$ we get the following $q$-analog of identity (\ref{identity4.1})
\begin{equation}\label{identity4.2}
\prod_{j=1}^l\bigg[{x+d_j\atop d_j}\bigg]_q
=\sum_{k=1}^dB_{\mathbf{d},k}(q)\bigg[{x+1\atop k}\bigg]_q.
\end{equation}
\begin{theorem}
A $q$-analog of identity (\ref{identity4.1}) is
$$\prod_{j=1}^l\bigg[{x+d_j\atop d_j}\bigg]_q
=\sum_{k=1}^dB_{\mathbf{d},k}(q)\bigg[{x+1\atop k}\bigg]_q.$$
\end{theorem}
Note that identity (\ref{identity4.2}) is a generalization of the following identity \cite{carlitz3}
$$\bigg[{x\atop 1}\bigg]_q^d=\sum_{k=0}^dq^{\frac{1}{2}k(k-1)}a_{d,k}(q)\bigg[{x\atop k}\bigg]_q$$
where $a_{d,k}(q)$ is a polynomial of $q$.

\section{Ordered Stirling numbers of the third kind over a multiset}

To derive a multiset version of identity (\ref{identity1.3}), we will compute $\sum_{\sigma\in\mathfrak{S}(\mathbf{d})}|\mathbb{Z}(n\alpha^d(\sigma))|$ by two different ways and to obtain a $q$-analog of this identity, we will compute $f_2(q)=\sum_{\sigma\in\mathfrak{S}(\mathbf{d})}\sum_{\mathbf{x}\in\mathbb{Z}(n\alpha^d(\sigma))}q^{\mathbf{x}}$ by two different ways.

Let $\mathbf{x}$ be a point in $\mathbb{Z}(n\alpha^{\mathbf{d}})$. Since $\alpha^{\mathbf{d}}$ is the disjoint union of the sets $I(\alpha^k)$ for $\alpha^k\in\mathcal{I}(\alpha^{\mathbf{d}})$, there is a unique simplex $\alpha^k$ in $\mathcal{I}(\alpha^{\mathbf{d}})$ such that $\mathbb{Z}(nI(\alpha^k))$ contains $\mathbf{x}$. If $\alpha^d(\sigma)$ is a $d$-dimensional simplex in $\mathcal{I}(\alpha^{\mathbf{d}})$ that contains $\alpha^k$, then we can represent $(\alpha^d(\sigma),\alpha^k)$ by a $k$-tuple of the form
\begin{equation}\label{identity5.1}
((\sigma_1,\ldots,\sigma_{i_1}),(\sigma_{i_1+1},\ldots,\sigma_{i_2}),\ldots,(\sigma_{i_{k-1}+1},\ldots,\sigma_d)),
\end{equation}
By a simple computation, the number of such $k$-tuples is ${d\choose d_1,d_2,\ldots,d_l}{d-1\choose k-1}$.
Note that the vertex set of $\alpha^k$ is
$\{\mathbf{e}_{i_{h,1},i_{h,2},\ldots,i_{h,l}}\,|\, h\in[k]_0\}$ where
$$\begin{cases}
(i_{0,1},i_{0,2},\ldots,i_{0,l})=\mathbf{0}\\
(i_{h,1},i_{h,2},\ldots,i_{h,l})=\sum_{i=1}^{i_h}\mathbf{e}_{l,i}\ \text{for}\ h\in[k-1]\\
(i_{k,1},i_{k,2},\ldots,i_{k,l})=(d_1,d_2,\ldots,d_l)
\end{cases}.$$
Therefore if we denote by $\big\lfloor{\mathbf{d}\atop k}\big\rfloor_O$ the number of $k$-tuples of form (\ref{identity5.1}), which is called an \emph{ordered Stirling number of the third kind} (or an \emph{ordered Lah number}) over a multiset, then we obtain
\begin{align*}
\sum_{\sigma\in\mathfrak{S}(\mathbf{d})}|\mathbb{Z}(n\alpha^d(\sigma))|
&=\sum_{\sigma\in\mathfrak{S}(\mathbf{d})}\sum_{k=1}^d\sum_{\alpha^k\subseteq
\alpha^d(\sigma)}
|\mathbb{Z}(nI(\alpha^k))|\\
&=\sum_{k=1}^d\bigg\lfloor{\mathbf{d}\atop k}\bigg\rfloor_O{n+1\choose k}.
 \end{align*}
In addition, 
$\sum_{\sigma\in\mathfrak{S}(\mathbf{d})}
|\mathbb{Z}(n\alpha^d(\sigma))|={d\choose d_1,d_2,\ldots,d_l}{n+d\choose d}$, thus it follows that
$${d\choose d_1,d_2,\ldots,d_l}{n+d\choose d}=\sum_{k=1}^d\bigg\lfloor{\mathbf{d}\atop k}\bigg\rfloor_O{n+1\choose k}.$$
As a result, we get a multiset version of identity (\ref{identity1.3})
\begin{equation}\label{identity5.1.1}
{d\choose d_1,d_2,\ldots,d_l}{x+d\choose d}=\sum_{k=1}^d\bigg\lfloor{\mathbf{d}\atop k}\bigg\rfloor_O{x+1\choose k}.\end{equation}
\begin{theorem}
Ordered Stirling numbers of the third kind over $S(\mathbf{d})$ satisfy
$$ {d\choose d_1,d_2,\ldots,d_l}{x+d\choose d}=\sum_{k=1}^d\bigg\lfloor{\mathbf{d}\atop k}\bigg\rfloor_O{x+1\choose k}.$$
\end{theorem}

Now we consider  $f_2(q)$.
The number of permutations of $S(\mathbf{d})$ is ${d\choose d_1,d_2,\ldots,d_l}$ and $\sum_{\mathbf{x}\in\mathbb{Z}(n\alpha^d(\sigma))}q^{\mathbf{x}}=\Big[{n+d\atop d}\Big]_q$, thus it follows that
\begin{equation}\label{identity5.2}
f_2(q)={d\choose d_1,d_2,\ldots,d_l}
\bigg[{n+d\atop d}\bigg]_q.
\end{equation}
In addition, by the second decomposition of $\alpha^{\mathbf{d}}$ we can obtain
\begin{align}\label{identity5.3}
f_2(q)=&\sum_{\sigma\in\mathfrak{S}(\mathbf{d})}\sum_{\mathbf{x}\in\mathbb{Z}(n\alpha^d(\sigma))}q^{\mathbf{x}}\nonumber\\
=&\sum_{\sigma\in\mathfrak{S}(\mathbf{d})}\sum_{k=1}^d\sum_{\alpha^k\in\mathcal{I}(\alpha^{\mathbf{d}})\atop
\alpha^k\subseteq\alpha^d(\sigma)}
\sum_{\mathbf{x}\in\mathbb{Z}(nI(\alpha^k))}q^{\mathbf{x}}\nonumber\\
=&\sum_{k=1}^d\sum_{\sigma\in\mathfrak{S}(\mathbf{d})}\sum_{\alpha^k\in\mathcal{I}(\alpha^{\mathbf{d}})\atop
\alpha^k\subseteq\alpha^d(\sigma)}
q^{maj(\alpha^k)}\bigg[{n+1\atop k}\bigg]_q\nonumber\\
=&\sum_{k=1}^dC_{\mathbf{d},k}(q)\bigg[{n+1\atop k}\bigg]_q
\end{align}
where $C_{\mathbf{d},k}(q)=\sum_{\sigma\in\mathfrak{S}(\mathbf{d})}\sum_{\alpha^k\in\mathcal{I}(\alpha^{\mathbf{d}})\atop
\alpha^k\subseteq\alpha^d(\sigma)}
q^{maj(\alpha^k)}$.
Therefore by combining identities (\ref{identity5.2}) and (\ref{identity5.3}) we obtain
$${d\choose d_1,d_2,\ldots,d_l}
\bigg[{n+d\atop d}\bigg]_q
=\sum_{k=1}^dC_{\mathbf{d},k}(q)\bigg[{n+1\atop k}\bigg]_q.$$
As a result, we get the following $q$-analog of identity (\ref{identity5.1.1})
$${d\choose d_1,d_2,\ldots,d_l}
\bigg[{x+d\atop d}\bigg]_q
=\sum_{k=1}^dC_{\mathbf{d},k}(q)\bigg[{x+1\atop k}\bigg]_q.$$
\begin{theorem}
A $q$-analog of identity (\ref{identity5.1.1}) is 
$${d\choose d_1,d_2,\ldots,d_l}
\bigg[{x+d\atop d}\bigg]_q
=\sum_{k=1}^dC_{\mathbf{d},k}(q)\bigg[{x+1\atop k}\bigg]_q.$$
\end{theorem}

\section*{Classification codes}
05A05, 05A10, 05A15, 05A18, 05A19, 05A30, 52B11


\begin{thebibliography}{10}

\bibitem{worpitzky} J. Worpitzky, Studien \"{u}ber die
Bernoullischen und Eulerschen Zaheln, Journal f.d. reine u. angew.
Math. 94 (1883) 203-232.

\bibitem{carlitz1} L. Carlitz, A combinatorial property of q-Eulerian numbers, Amer. Math. Monthly 82 No. 1 (1975) 51-54.

\bibitem{carlitz2} L. Carlitz, $q$-Bernoulli and Eulerian numbers.
Trans. Amer. Math. Soc. 76 (1954) 332-350.

\bibitem{tsylova} E. G. Tsylova, The asymptotic behavior of generalized Stirling numbers, Combinatorial-algebraic methods in applied math, Gorkov. Gos. Univ., Gorki, 143-154, 158, 1985.


\bibitem{aigner} M. Aigner, Combinatorial theory, Springer, New York, 1979.

\bibitem{motzkin} T. Motzkin, Sorting numbers for cylinders and other classification numbers, Proc. Symp. Pure Math., Vol. 19 (1971) 167-176.

\bibitem{lah} I. Lah, Eine neue Art von Zahlen, ihre Eigenschaften und Anwendung in der mathematischen Statistik, Mitteilungsbl. Math. Statist. 7 (1955) 203-212.

\bibitem{stanley} R. P. Stanley, Eulerian partitions of a unit hypercube, Higher Combinatorics (M. Aigner, ed.), Reidel, Dordrecht/Boston, 1977, pp. 49. 


\bibitem{andrews} G. E. Andrews, The theory of partitions,
Cambridge University Press, Cambridge, 1998.

\bibitem{sylvester} J. J. Sylvester, A constructive theory of partitions, arranged in three acts, an interact and an exodian, Amer. J. Math. 5 (1882) 251-330.

\bibitem{carlitz3}  L. Carlitz, $q$-Bernoulli numbers and polynomials. Duke Math. J. 15 (1948) 987–1000.

\end{thebibliography}
\end{document}